# Taking all positive eigenvectors is suboptimal in classical multidimensional scaling

Jeffrey Tsang, Rajesh Pereira

jeffrey.tsang@ieee.org  pereirar@uoguelph.ca

## Abstract

It is hard to overstate the importance of multidimensional scaling as an analysis technique in the broad sciences. Classical, or Torgerson multidimensional scaling is one of the main variants, with the advantage that it has a closed-form analytic solution. However, this solution is exact if and only if the distances are Euclidean. Conversely, there has been comparatively little discussion on what to do in the presence of negative eigenvalues: the intuitive solution, prima facie justifiable in least-squares terms, is to take every positive eigenvector as a dimension. We show that this, minimizing least-squares to the centred distances instead of the true distances, is suboptimal — throwing away positive eigenvectors can decrease the error even as we project to fewer dimensions. We provide provably better methods for handling this common case.

## 1 Introduction

Multidimensional scaling is a fundamental analysis technique that takes as input a matrix of distances or dissimilarities between items and returns a configuration of points in Euclidean space such that the inter-point distances approximate the input. From the original definition and usage in psychology [12, 22], applications have spread throughout all the sciences, for example in ecology [10], genetics [13], palynology [17], neuroscience [7], medicine [3], education [21], management [19], music [8], physical chemistry [25], and electrical engineering [9]. Multiple textbooks have been written on the subject, see for example [1, 4].

Classical multidimensional scaling was the first version, which operates based on matrix eigendecomposition. This gives it the advantage of having a closed-form solution; in the case that no exact solution exists, an intuitive solution with a seemingly straightforward error analysis can be found. Unfortunately, that analysis is heavily marred by the simple fact that it relies on minimizing least-squares to the centred distances instead of the true distances. We show that the naïve solution is therefore suboptimal and provide better methods for handling this common case.

## 2   Preliminaries

Unless otherwise noted, all matrices will be of size $n \times n$. We will freely switch between capital letters and entry-form $[a_{ij}]$, that is the $i,j$th entry of the matrix is $a_{ij}$, as a function of $i$ and $j$. Let $\delta(i = j)$ be the Kronecker delta function, which is 1 if $i = j$ and 0 otherwise. Let $I = [\delta(i = j)]$ denote the identity matrix, $J = [1]$ denote the matrix of all ones, and $\Phi$ be the linear operator that takes a matrix and zeros out all entries off the main diagonal. That is,

$$\Phi\left([a_{ij}]\right) = [\delta(i = j)a_{ij}].$$

We will follow the exposition in [11]. Let $D = [d_{ij}^2]$ be the input data matrix of squared distances: we will presume that it is symmetric, has a zero diagonal and otherwise nonnegative entries, and the (non-squared) entries satisfy the triangle inequality. To prepare for classical MDS, we will double-centre the matrix and scale this by $-1/2$; call this $B$:

$$B = [b_{ij}] = \left[-\frac{1}{2}(d_{ij}^2 - d_{i\cdot}^2 - d_{\cdot j}^2 + d_{\cdot\cdot}^2)\right] = -\frac{1}{2}\left(I - \frac{1}{n}J\right) D \left(I - \frac{1}{n}J\right),$$

where $d_{i\cdot}^2$ is the mean of the $i$th row of the matrix, $d_{\cdot j}^2$ the mean of the $j$th column, and $d_{\cdot\cdot}^2$ the mean of all $n^2$ entries of the matrix.

To see the point of this transformation, assume for now that the matrix comes from a set of $n$ vectors in $\mathbb{R}^n$, say $\{\vec{x}_1, \vec{x}_2, \ldots, \vec{x}_n\}$. Since distances are invariant under translation, further suppose the set of vectors has mean $\bar{x} = \frac{1}{n}\sum_{k=1}^n \vec{x}_k = \vec{0}$. So $d_{ij}^2 = \|\vec{x}_i - \vec{x}_j\|^2 = \|\vec{x}_i\|^2 - 2\langle \vec{x}_i, \vec{x}_j\rangle + \|\vec{x}_j\|^2$. Thus,

$$-2b_{ij} = \|\vec{x}_i - \vec{x}_j\|^2 - \frac{1}{n}\sum_{k=1}^n \|\vec{x}_i - \vec{x}_k\|^2 - \frac{1}{n}\sum_{k=1}^n \|\vec{x}_j - \vec{x}_k\|^2 + \frac{1}{n^2}\sum_{k,l=1}^{n,n} \|\vec{x}_k - \vec{x}_l\|^2$$

$$= \left(\|\vec{x}_i\|^2 - 2\langle \vec{x}_i, \vec{x}_j\rangle + \|\vec{x}_j\|^2\right) - \left(\|\vec{x}_i\|^2 - 2\langle \vec{x}_i, \bar{x}\rangle + \frac{1}{n}\sum_{k=1}^n \|\vec{x}_k\|^2\right)$$

$$- \left(\|\vec{x}_j\|^2 - 2\langle \vec{x}_j, \bar{x}\rangle + \frac{1}{n}\sum_{k=1}^n \|\vec{x}_k\|^2\right) + \left(\frac{2}{n}\sum_{k=1}^n \|\vec{x}_k\|^2 - 2\langle \bar{x}, \bar{x}\rangle\right)$$

$$= -2\langle \vec{x}_i, \vec{x}_j\rangle.$$

Since $b_{ij} = \langle \vec{x}_i, \vec{x}_j\rangle$, $B$ is a Gram matrix of inner products, formed by $XX^T$, where $X = [\vec{x}_1, \vec{x}_2, \ldots, \vec{x}_n]$, the matrix formed by concatenating each $\vec{x}_i$ as a column.

We now attempt to take the square root of the matrix to recover the vectors. To do so, we diagonalize the matrix:

$$B = V\Lambda V^T = \left(V\Lambda^{\frac{1}{2}}\right)\left(\Lambda^{\frac{1}{2}}V^T\right) = \left(V\Lambda^{\frac{1}{2}}\right)\left(V\Lambda^{\frac{1}{2}}\right)^T$$

where the eigenvectors and singular vectors coincide as the matrix is symmetric. To have a real solution, we therefore need that all eigenvalues of $B$ be nonnegative. Let us list the eigenvalues of $B$ in descending order as $\lambda_1 \geq \lambda_2 \geq \ldots \geq \lambda_n$.

Note that we could write a linear transformation to turn $D$ into $B$; we will write out its inverse here.

$$b_{ii}+b_{jj}-2b_{ij} = (d_{ii}^2-d_{i\cdot}^2-d_{\cdot i}^2+d_{\cdot\cdot}^2)+(d_{jj}^2-d_{j\cdot}^2-d_{\cdot j}^2+d_{\cdot\cdot}^2)-2(d_{ij}^2-d_{i\cdot}^2-d_{\cdot j}^2+d_{\cdot\cdot}^2) = d_{ij}^2$$

since by symmetry $d_{i\cdot}^2 = d_{\cdot i}^2$. We have that $\Phi(B)J = [\sum_{k=1}^{n}(\delta(i=k)b_{ik})(1)] = [b_{ii}]$ and similarly that $J\Phi(B) = [b_{jj}]$. Thus the equation

$$D = \Phi(B)J + J\Phi(B) - 2B,$$

completing the preliminaries.

## 3 Least-squares error

Consider the problem of having negative eigenvalues in $B$. A celebrated theorem [20] proves that there cannot be any exact solution to MDS in this case, so we need to define a measure of the error to optimize. Naturally, some form of least-squares is indicated, for example

$$L(X) = \sum_{i,j=1}^{n,n} \left(d_{ij}^2 - \|\vec{x_i} - \vec{x_j}\|^2\right)^2$$

on the squared distances, which defines classical MDS. We will not go into detail regarding the competing variants (for example, defining least-squares on the non-squared distances leads to a different objective, see [5, 6]).

From now on, let $X$ be an approximate solution to MDS. For now, consider the ur-error

$$L'(X) = \sum_{i,j=1}^{n,n} \left(b_{ij} - (XX^T)_{ij}\right)^2$$

which is the sum of the squares of each entry of $B - XX^T$, also known as the square of the Frobenius norm, or $\left\|B - XX^T\right\|_F^2 = \text{tr}((B-XX^T)(B-XX^T)^T)$. The trace of a matrix is the sum of its diagonal, $\text{tr}([a_{ij}]) = \sum_{i=1}^{n} a_{ii}$. Standard results in matrix analysis state that the Frobenius norm squared is the sum of squares of the singular values (= sum of squares of eigenvalues for a symmetric matrix).

Thus the intuitive solution to the negative eigenvalue problem above: pick $XX^T = V\Lambda_+V^T$, that is, let $XX^T$ and $B$ have the same eigenvectors, with all negative eigenvalues clamped to 0 [14]. The minimized "error" in this case is simply the sum of squares of negative eigenvalues.

Unfortunately, the error is phrased in terms of the entries of $D$, not $B$. Since $D = \Phi(B)J + J\Phi(B) - 2B$, in the same way, the approximated distances are given by $\Phi(XX^T)J + J\Phi(XX^T) - 2XX^T$. Thus the true error is

$$L(X) = \left\|(\Phi(XX^T)J + J\Phi(XX^T) - 2XX^T) - (\Phi(B)J + J\Phi(B) - 2B)\right\|_F^2$$
$$= \left\|(\Phi(B - XX^T)J + J\Phi(B - XX^T)) - 2(B - XX^T)\right\|_F^2$$

by linearity of $\Phi$. Write as shorthand $R = [r_{ij}] = B - XX^T$.

$$= \|\Phi(R)J + J\Phi(R)\|_F^2 - 2\operatorname{tr}\left((2R)(\Phi(R)J + J\Phi(R))^T\right) + \|2R\|_F^2$$

from the trace definition, where trace is invariant under transposition. Furthermore, it is invariant under switching the order, so

$$= \|\Phi(R)J + J\Phi(R)\|_F^2 - 4\operatorname{tr}\left(RJ\Phi(R) + \Phi(R)JR\right) + 4\|R\|_F^2$$

Note that both $B$ and $XX^T$ by extension have zero column and row sums, due to the centering process. Since multiplying by $J$ is summing rows or columns, $JR = RJ = [0]$. Thus the cross-term cancels:

$$= \|\Phi(R)J + J\Phi(R)\|_F^2 + 4\|R\|_F^2.$$

The first term in the above can now be written as

$$\sum_{i,j=1}^{n,n} (r_{ii} + r_{jj})^2 = n\sum_{i=1}^{n} r_{ii}^2 + 2\left(\sum_{i=1}^{n} r_{ii}\right)\left(\sum_{j=1}^{n} r_{jj}\right) + n\sum_{j=1}^{n} r_{jj}^2$$

$$= 2n\sum_{i=1}^{n} r_{ii}^2 + 2(\operatorname{tr}(R))^2.$$

With that, the error simplifies to

$$L(X) = 2n\sum_{i=1}^{n} r_{ii}^2 + 2(\operatorname{tr}(R))^2 + 4\|R\|_F^2.$$

## 4 Error bounds

Although the first term cannot be easily simplfied, we can bound it. First note that we can find $\operatorname{tr}(R)$ by summing the diagonal entries. To minimize the sum of squares, we can set each diagonal entry equal, to $\operatorname{tr}(R)/n$, whence the sum is bounded below by $n \times (\operatorname{tr}(R))^2/n^2$, and the first term by $2(\operatorname{tr}(R))^2$.

To maximize, consider that $\|R\|_F^2$ is the sum of squares of all entries, clearly not less than the sum of squares of the diagonal. We also know that $R$ has row and column sums zero. Thus to maximize the ratio between $\sum_{i=1}^{n} r_{ii}^2$ and $\|R\|_F^2$, let us set each diagonal entry equal (to say $c$), and each off-diagonal entry equal (to $-\frac{c}{n-1}$). The sum we wish to bound is $nc^2$, and the upper bound $nc^2 + n(n-1)\frac{c^2}{(n-1)^2} = (nc^2)(1 + \frac{1}{n-1})$, a factor of $\frac{n}{n-1}$. Therefore the first term is bounded above by $2(n-1)\|R\|_F^2$. In conclusion, we have

$$4(\operatorname{tr}(R))^2 + 4\|R\|_F^2 \leq L(X) \leq 2(\operatorname{tr}(R))^2 + (2 + 2n)\|R\|_F^2.$$

Our bounds are tight, since the equal-diagonal, equal-off-diagonal matrix is realizable by the regular simplex, that is, set all distances to 1 ($D = J - I$).

After centering, ignore $X$ (set every MDS point to $\vec{0}$) and let $B = R$ — then the lower and upper bounds coincide.

Note that $(\text{tr}(R))^2$ is not the same as $\|R\|_F^2$ — the first is the square of the sum of eigenvalues, the second is the sum of the squares. In fact, let us analyze from that viewpoint. The latter term is the ur-error minimized by $V\Lambda_+ V^T$. The former term squares the sum of residual eigenvalues, which means to minimize it we have to leave enough positive eigenvalues to balance out the sum (as we cannot duplicate negative eigenvalues).

Therefore, to minimize a weighted sum of both terms, we cannot simply match every single positive eigenvalue. However, we would be seriously remiss if we failed to mention that even if a function can be bounded both above and below, the minimum of the function need not be between the minima of the upper and lower bounds. An example is depicted in Figure 1.

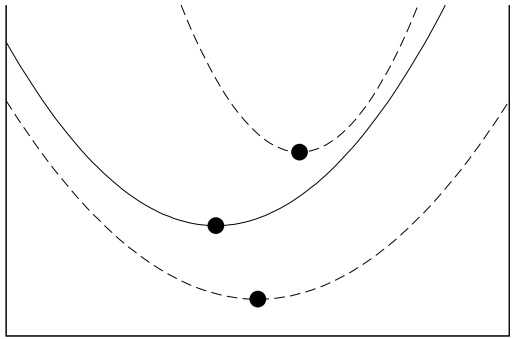

**Figure 1:** Graphical demonstration that $f(x) \leq g(x) \leq h(x)\ \forall x$ does not imply $\operatorname{argmin}_x f(x) \leq \operatorname{argmin}_x g(x) \leq \operatorname{argmin}_x h(x)$.

## 5 Improved algorithms

Leaving that caveat aside for now, we can consider simple algorithms to minimize these bounds. Let us treat the problem as two-dimensional: we can trivially compute $(\text{tr}(R))^2$ and $\|R\|_F^2$, in terms of the eigenvalues of $X = V\Lambda'^{\frac{1}{2}}$, where $V = [v_{ij}]$ is the matrix of eigenvectors of $B$ and $\Lambda'$ is the diagonal matrix of reconstructed eigenvalues. Let its entries be $\lambda_1', \lambda_2', \ldots, \lambda_n' \geq 0$. Then $(\text{tr}(R))^2 = \left(\sum_{i=1}^n \lambda_i - \lambda_i'\right)^2$ and $\|R\|_F^2 = \sum_{i=1}^n (\lambda_i - \lambda_i')^2$, which is entirely in terms of the $\lambda_i'$s.

We solve this with a marginalization procedure. Suppose we fix $(\text{tr}(R))^2$, or equivalently fix $\sum_{i=1}^n \lambda_i'$. How would we minimize $\|R\|_F^2$? Remember that we are summing up the squares of the unmatched eigenvalues: thus the minimum is when we equalize the maximal unmatched terms. That is, fix some $c \in \mathbb{R}$,

and pick
$$\lambda'_i = \begin{cases} \lambda_i - c, & \lambda_i \geq c \\ 0, & \lambda_i < c \end{cases}.$$

This solution essentially caps all eigenvalues of $R$ at $c$, leaving lesser ones unchanged — by monotonicity there exists a unique $c$. Now note that the bounds for $L(X)$ can be written as a piecewise quadratic function in $c$ (as $\text{tr}(R)$ changes linearly in $c$ between different $\lambda_i$), and this can be easily solved in each piece.

Given a candidate solution $X$, we can compute the true value of $L(X)$ in quadratic time: testing the $\mathcal{O}(n)$ candidates and picking the best one takes cubic time, which is the same complexity as diagonalizing $B$ in the first place. Hence the algorithm takes no longer asymptotically than the usual one. Furthermore, as the space of solutions considered includes the original ($c = 0$ matches all positive eigenvalues), this solution clearly dominates the original.

We outline a second algorithm that avoids the complications illustrated in Figure 1. We can directly write out $\sum_{i=1}^n r_{ii}^2$ as a function of the $\lambda'_i$s:

$$r_{ii} = (B - V\Lambda'V^T)_{ii} = b_{ii} - \sum_{k=1}^n \lambda'_k v_{ik}^2$$

where $b_{ii}$, $v_{ik}^2$ are constants found from diagonalizing $B$. Hence $\sum_{i=1}^n r_{ii}^2$ is a quadratic function of all the $\lambda'_i$s, along with $(\text{tr}(R))^2$ and $\|R\|_F^2$ which we have treated above. Rephrasing, we wish to minimize $L(X(\lambda'_1, \lambda'_2, \ldots, \lambda'_n))$, a quadratic function, under the constraint that $\lambda'_i \geq 0$.

This is a quadratic program. Not only that, it is relatively easy to show that $L(X)$, as a sum of squares, will be a convex quadratic program, whence a cubic algorithm exists to solve it [16]. Again, the space of feasible solutions includes the original, and in cubic (that is, asymptotically equal) time, we have an algorithm that solves classical MDS and dominates the original in least-squares error.

To facilitate implementation, let us directly expand $L(X)$. Typically, the function is written in the form $\frac{1}{2}\vec{x}^T Q \vec{x} + \vec{c} \cdot \vec{x} + a$, where $Q$ is symmetric, and the constant $a$ is immaterial to the optimization and usually left out.

$$L(X) = 2n \sum_{i=1}^n r_{ii}^2 + 2(\text{tr}(R))^2 + 4\|R\|_F^2$$
$$= 2n \sum_{i=1}^n \left(b_{ii} - \sum_{k=1}^n \lambda'_k v_{ik}^2\right)^2 + 2\left(\sum_{i=1}^n \lambda_i - \sum_{i=1}^n \lambda'_i\right)^2 + 4\sum_{i=1}^n (\lambda_i - \lambda'_i)^2$$

$$
\begin{aligned}
&= 2n \sum_{i,j=1}^{n,n} \left( \sum_{k=1}^{n} v_{ki}^2 v_{kj}^2 \right) \lambda_i' \lambda_j' - 4n \sum_{i=1}^{n} \left( \sum_{k=1}^{n} b_{kk} v_{ki}^2 \right) \lambda_i' + 2n \sum_{i=1}^{n} b_{ii}^2 \\
&\quad + 2 \left( \sum_{i=1}^{n} \lambda_i' \right)^2 - 4 \sum_{i=1}^{n} \left( \sum_{k=1}^{n} \lambda_k \right) \lambda_i' + 2 \left( \sum_{i=1}^{n} \lambda_i \right)^2 + 4 \sum_{i=1}^{n} \lambda_i'^2 \\
&\quad - 8 \sum_{i=1}^{n} \lambda_i \lambda_i' + 4 \sum_{i=1}^{n} \lambda_i^2 \\
&= 2n(\vec{\lambda}')^T W^T W \vec{\lambda}' - 4n W^T \operatorname{diag}(B) \cdot \vec{\lambda}' + 2n \sum_{i=1}^{n} b_{ii}^2 + 2(\vec{\lambda}')^T J \vec{\lambda}' \\
&\quad - 4 \operatorname{tr}(B) \vec{1} \cdot \vec{\lambda}' + 2 \left( \operatorname{tr}(B) \right)^2 + 4(\vec{\lambda}')^T I \vec{\lambda}' - 8 \operatorname{diag}(\Lambda) \cdot \vec{\lambda}' + 4 \| B \|_F^2
\end{aligned}
$$

where $W = [v_{ij}^2]$ is the entry-wise square of the eigenvectors matrix $V$, $\operatorname{diag}(B)$ is the vectorized diagonal of $B$, $\vec{1}$ is the vector of all ones, $\operatorname{diag}(\Lambda)$ is the vector of eigenvalues of $B$ in order, and $\vec{\lambda}'$ are our reconstructed eigenvalues (the parameters of solution). Therefore, we have

$$
\begin{aligned}
Q &= 4n W^T W + 4J + 8I \\
\vec{c} &= -4n W^T \operatorname{diag}(B) - 4 \operatorname{tr}(B) \vec{1} - 8 \operatorname{diag}(\Lambda) \\
a &= 2n \sum_{i=1}^{n} b_{ii}^2 + 2 \left( \operatorname{tr}(B) \right)^2 + 4 \| B \|_F^2
\end{aligned}
$$

and after solving for $\vec{\lambda}'$, our $\Lambda'$ is the diagonal matrix of $\vec{\lambda}'$, and the MDS projected points are given by $V \Lambda'^{\frac{1}{2}}$.

The only remaining issue is that this presumes the solution should have the same eigenvectors as $B$, which is no longer necessarily true as it was for $L'(X)$. We leave this as an open question for future work.

## 6 Example calculations

We provide an example from real-world usage to demonstrate that the preceding is no mere theoretical problem. We have a distance matrix of size $296 \times 296$, as computed following [23], using the updated algorithm in [24]. Other than taking the dataset and noting that it is explicitly known to be non-Euclidean, we bear no connection to that work.

We go through the motions of running classical MDS on this dataset, and tabulate a variety of statistics on the output in Table 1. Aside from the eigenvalues in decreasing order, we will show the normalized least-squares error, which can be directly interpreted as the root mean square error on the distances themselves. For a sense of scale, the input distances are bounded above by 1, and have a mean of 0.288577 and RMS of 0.323568.

All following algorithms are forced to consider only the top $n$ eigenvectors ($\lambda_i' \equiv 0$ for $i > n$): the standard algorithm copies the top $n$ positive eigenvalues;

we can optimize the upper and lower bounds in terms of $c$, the eigenvalue cutoff, then take that as a solution; we also present the quadratic programming solution. For comparison, we also present the best SMACOF solution in $n$ dimensions over 100 random initial configurations. SMACOF is an iterative algorithm that minimizes the least-squares error on the non-squared distances [5].

The first major insight is that even though there are 103 positive eigenvalues, using any more than just the top 5 eigenvalues starts *increasing* the error even as more dimensions are used. In fact, using all 103 is worse than simply using 3! This extremely small number, compared to the number of positive eigenvalues, shows that the problem is of clear practical relevance.

Note that both the upper and lower bounds are uniformly better than the standard solution: the cutoff starts off negative, as when forced to use insufficient dimensions, it is gainful to intentionally overscale the eigenvectors to increase the average distances. The upper bound fares much worse as the $2+2n$ factor on the $\|R\|_F^2$ term means it is hardly different from the original problem and thus its error eventually boomerangs; however it uses only the top 49 eigenvectors. The lower bound improves the best solution by 14% and uses just 10 components.

The quadratic programming solution beats them all, improving by 18%. Peculiarly, it refuses to use the 11th, 13th, and 15–18th eigenvectors, and none of the positive ones after the 20th. Even more weirdly, it somehow uses the 246th and 258th eigenvectors and only those, even though their eigenvalues are negative! This is not a numerical artifact: their reconstructed eigenvalues are 0.003882 and 0.006092; the error is also clearly decreasing.

The SMACOF algorithm technically optimizes a different objective, and hence the comparison in $L(X)$ terms is not completely fair. Even then, it manages to outperform the quadratic program in 2–4 dimensions. When more dimensions are used, the combinatorial explosion causes the algorithm to get stuck in local minima and improvement eventually vanishes. We do not attempt to run it for $n > 21$ dimensions.

Thus under the standard algorithm, the data seems to be 5-dimensional; using the lower bound, an improved fit is found with 10 dimensions, and with the optimal quadratic programming method, 16 dimensions are used. These are maximal fits, in the sense that even if allowed extra dimensions, the methods will not use them.

To close, we note that two of the most widely used statistical software packages are unaware of this problem. MATLAB [15] acknowledges that in the presence of significant negative eigenvalues, its choice of using all positive eigenvectors may yield a poor solution. R [18] provides several options, including using a minimal additive constant from [2] to make the distance matrix Euclidean. However, it still follows the work of [14] and hence also uses all positive eigenvectors by default.

## Supplementary Materials:

D2ST.dist: a non-squared distance matrix of size $296 \times 296$, given in upper triangular order, listing $d_{1,1}, d_{1,2}, \ldots, d_{1,296}, d_{2,3}, \ldots$ (total 43660 entries); the numbers are stored in IEEE 754 double precision format. Available from arXiv→download→source.

## Acknowledgements

Rajesh Pereira acknowledges the support of a Natural Sciences and Engineering Council of Canada Discovery Grant.

|   | | Root mean square error on $d_{ij}^2$, $100 \times \sqrt{L(X)/n(n-1)}$ | | | | |
|---|---|---|---|---|---|---|
| # | Eigenvalue | Standard | Upper(Cutoff) | Lower(Cutoff) | CQP | SMACOF |
| 1 | 12.870775 | 2.20519329 | 2.20178706(−0.00863) | 2.14970170(−1.28598) | 2.05734225 | 2.21289091 |
| 2 | 1.037833 | 1.45592702 | 1.45178684(−0.00513) | 1.31937723(−0.51137) | 1.29080402 | 1.26708425 |
| 3 | 0.841390 | 0.89622191 | 0.89409993(−0.00231) | 0.83180580(−0.17318) | 0.82418493 | 0.78480269 |
| 4 | 0.476908 | 0.65701404 | 0.65680429(−0.00072) | 0.65854993(−0.04316) | 0.65147324 | 0.62506041 |
| 5 | 0.383134 | 0.55749502 | 0.55692006(0.000554) | 0.53853137(0.027886) | 0.52684668 | 0.53871108 |
| 6 | 0.240600 | 0.56075204 | 0.55774720(0.001346) | 0.48723979(0.058274) | 0.47279706 | 0.49638691 |
| 7 | 0.124567 | 0.60165224 | 0.59623265(0.001751) | 0.48338972(0.066560) | 0.46840905 | 0.49234572 |
| 8 | 0.095117 | 0.63479317 | 0.62693049(0.002058) | 0.48102608(0.069733) | 0.46517676 | 0.49090814 |
| 9 | 0.090882 | 0.67155736 | 0.66077926(0.002348) | 0.47911194(0.071848) | 0.46186658 | 0.49012857 |
| 10 | 0.083647 | 0.70978791 | 0.69587863(0.002613) | 0.47876356(0.072921) | 0.46080720 | 0.48968884 |
| 11 | 0.066470 | 0.74715742 | 0.72995180(0.002820) | (same) | (same) | 0.48925221 |
| 12 | 0.060986 | 0.77286918 | 0.75259874(0.003008) | (same) | 0.45690971 | 0.48943608 |
| 13 | 0.054237 | 0.80406322 | 0.78039058(0.003174) | (same) | (same) | 0.48955872 |
| 14 | 0.045480 | 0.82707314 | 0.80020239(0.003310) | (same) | 0.45622298 | 0.48961516 |
| 15 | 0.042569 | 0.85124591 | 0.82100872(0.003435) | (same) | (same) | 0.48968960 |
| 16 | 0.039933 | 0.87451480 | 0.84080891(0.003552) | (same) | (same) | 0.48972780 |
| 17 | 0.036081 | 0.89475692 | 0.85769606(0.003656) | (same) | (same) | 0.48975081 |
| 18 | 0.035816 | 0.91614215 | 0.87538331(0.003758) | (same) | (same) | 0.48973887 |
| 19 | 0.028725 | 0.93216642 | 0.88808096(0.003837) | (same) | 0.45594512 | 0.48978879 |
| 20 | 0.025852 | 0.94636821 | 0.89905662(0.003906) | (same) | 0.45580692 | 0.48982703 |

**Table 1:** Root mean square results for classical MDS performed using various methods, restricted to using only the top $n$ eigenvectors. Standard: copy all eigenvalues; Upper: minimize the upper bound; Cutoff: the $c$ value used; Lower: minimize the lower bound; CQP: optimize the least-squares quadratic program. SMACOF: use the iterative stress majorization algorithm in $n$ dimensions. (same) means the solution is exactly the same as for $n-1$ dimensions.

|     |            | Root mean square error on $d_{ij}^2$, $100 \times \sqrt{L(X)/n(n-1)}$ | | | | |
| --- | --- | --- | --- | --- | --- | --- |
| #   | Eigenvalue | Standard   | Upper(Cutoff)         | Lower(Cutoff) | CQP        | SMACOF     |
| 21  | 0.024536   | 0.96135228 | 0.91063699(0.003971)  | (same)        | (same)     | 0.48983074 |
| ⋮   | ⋮          | ⋮          | ⋮                     | ⋮             | ⋮          | ⋮          |
| 49  | 0.004537   | 1.14058061 | 1.00324282(0.004497)  | (same)        | (same)     | N/A        |
| 50  | 0.004050   | 1.14299177 | (same)                | (same)        | (same)     | N/A        |
| ⋮   | ⋮          | ⋮          | ⋮                     | ⋮             | ⋮          | ⋮          |
| 103 | 0.000018   | 1.18422460 | (same)                | (same)        | (same)     | N/A        |
| 104 | 0.000000   | (same)     | (same)                | (same)        | (same)     | N/A        |
| 105 | −0.000035  | (same)     | (same)                | (same)        | (same)     | N/A        |
| ⋮   | ⋮          | ⋮          | ⋮                     | ⋮             | ⋮          | ⋮          |
| 245 | −0.004660  | (same)     | (same)                | (same)        | (same)     | N/A        |
| 246 | −0.004816  | (same)     | (same)                | (same)        | 0.45578827 | N/A        |
| 247 | −0.004953  | (same)     | (same)                | (same)        | (same)     | N/A        |
| ⋮   | ⋮          | ⋮          | ⋮                     | ⋮             | ⋮          | ⋮          |
| 257 | −0.006288  | (same)     | (same)                | (same)        | (same)     | N/A        |
| 258 | −0.006408  | (same)     | (same)                | (same)        | 0.45572899 | N/A        |
| 259 | −0.006525  | (same)     | (same)                | (same)        | (same)     | N/A        |
| ⋮   | ⋮          | ⋮          | ⋮                     | ⋮             | ⋮          | ⋮          |
| 296 | −0.387203  | (same)     | (same)                | (same)        | (same)     | N/A        |

**Table 1** (continued): Root mean square results for classical MDS performed using various methods, restricted to using only the top $n$ eigenvectors. Standard: copy all eigenvalues; Upper: minimize the upper bound; Cutoff: the $c$ value used; Lower: minimize the lower bound; CQP: optimize the least-squares quadratic program. SMACOF: use the iterative stress majorization algorithm in $n$ dimensions. (same) means the solution is exactly the same as for $n-1$ dimensions.